# Computational Potential Energy Minimization Studies on the Prion AGAAAAGA Amyloid Fibril Molecular Structures

*Jiapu Zhang*

## 1. INTRODUCTION

X-ray crystallography, NMR (Nuclear Magnetic Resonance) spectroscopy, and dual polarization interferometry, etc are indeed very powerful tools to determine the 3D structure of a protein (including the membrane protein), though they are time-consuming and costly. However, for some proteins, due to their unstable, noncrystalline and insoluble nature, these tools cannot work. Under this condition, mathematical and physical theoretical methods and computational approaches allow us to obtain a description of the protein 3D structure at a submicroscopic level. This Chapter presents some practical and useful mathematical optimization computational approaches to produce 3D structures of the Prion AGAAAAGA Amyloid Fibrils, from an energy minimization point of view.

X-ray crystallography finds the X-ray final structure of a protein, which usually need refinements in order to produce a better structure. The computational methods presented in this Chapter can be also acted as a tool for the refinements.

All neurodegenerative diseases including Parkinson's, Alzheimer's, Huntington's, and Prion's have a similarity, which is they all featured amyloid fibrils (en.wikipedia.org/wiki/Amyloid and references (Nelson et al., 2005; Sawaya et al., 2007; Sunde et al., 1997; Wormell, 1954; Gilead and Gazit, 2004; Morley et al. 2006; Gazit, 2002; Pawar et al., 2005; and references therein). A prion is a misshapen protein that acts like an infectious agent (hence the name, which comes from the words protein and infection). Prions cause a number of fatal diseases such as mad cow disease in cattle, scrapie in sheep and kuru and Creutzfeldt-Jakob disease (CJD) in humans. Prion diseases (being rich in β-sheets (about 43% β-sheet) (Griffith, 1967; Cappaia and Collins, 2004; Daude, 2004; Ogayar and Snchez-Prez, 1998; Pan et al., 1993; Reilly, 2000) belong to neurodegenerative diseases. Many experimental studies such as (Brown, 2000; Brown, 2001; Brown, 1994; Cappai and Collins, 2004; Harrison et al., 2010; Holscher, 1998; Jobling et al., 2001; Jobling et al., 1999; Kuwata et al., 2003; Norstrom and Mastrianni, 2005; Wegner et al., 2002; Laganowsky et al., 2012; Jones et al., 2012; Sasaki et al., 2008; Haigh et al., 2005; Kourie et al., 2003; Zanuy et al., 2003; Kourie, 2001; Chabry et al., 1998; Gasset et al., 1992) have shown that the normal hydrophobic region (113-120) AGAAAAGA of prion proteins is an inhibitor/blocker of prion diseases. PrP lacking this palindrome could not convert to prion diseases. The presence of residues 119 and 120

(the two last residues within the motif AGAAAAGA) seems to be crucial for this inhibitory effect. The replacement of Glycine at residues 114 and 119 by Alanine led to the inability of the peptide to build fibrils but it nevertheless increased. The A117V variant is linked to the GSS disease. The physiological conditions such as pH (Cappai and Collins, 2004) and temperature (Wagoner et al., 2011) will affect the propensity to form fibrils in this region. The 3D atomic resolution structure of PrP (106-126), i.e. TNVKHVAGAAAAGAVVGGLGG, can be looked as the structure of a control peptide (Cheng et al., 2011; Lee et al., 2008). Ma and Nussinov (2002) established homology structure of AGAAAAGA and its molecular dynamics simulation studies. Recently, Wagoner et al. studied the structure of GAVAAAAVAG of mouse prion protein (Wagoner, 2010; Wagoner et al., 2011). Furthermore, the author computationally clarified that prion AGAAAAGA segment indeed has an amyloid fibril forming property (Fig. 1).

*Fig. 1: Prion AGAAAAGA (113-120) is surely and clearly identified as the amyloid fibril formation region, because its energy is less than the amyloid fibril formation threshold energy of -26 kcal/mol (Zhang et al., 2007).*

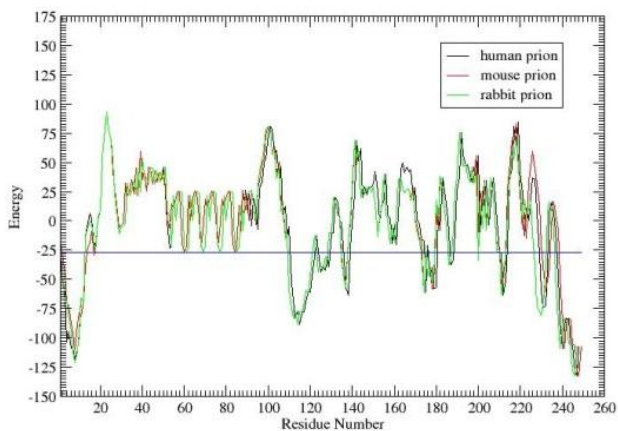

However, to the best of the author's knowledge, there is little X-ray or NMR structural data available to date on AGAAAAGA, which falls just within the N-terminal unstructured region (1–123) of prion proteins, due to its unstable, noncrystalline and insoluble nature. This Chapter will computationally study the molecular modeling (MM) structures of this region of prions.

## 2. MOLECULAR STRUCTURES OF PRION AGAAAAGA AMYLOID FIBRILS

"Amyloid is characterized by a cross-β sheet quaternary structure" and "recent X-ray diffraction studies of microcrystals revealed atomistic details of core region of amyloid" (en.wikipedia.org/wiki/Amyloid and references (Nelson et al., 2005; Sawaya et al., 2007; Sunde et al., 1997; Wormell, 1954; Gilead and Gazit, 2004; Morley et al., 2006; Gazit,

2002; Pawar et al., 2005; and references therein). All the quaternary structures of amyloid cross-β spines can be reduced to the one of 8 classes of steric zippers of (Sawaya et al., 2007), with strong van der Waals (vdw) interactions between β-sheets and hydrogen bonds (HBs) to maintain the β-strands.

A new era in the structural analysis of amyloids started from the 'steric zipper'- β-sheets (Nelson et al., 2005). As the two sheets zip up, HPs (Hydrophobic Packings) (& vdws) have been formed. The extension of the 'steric zipper' above and below (i.e. the β-strands) is maintained by HBs (but there is no HB between the two β-sheets). This is the common structure associated with some 20 neurodegenerative amyloid diseases, ranging from Alzheimer's and type-II diabetes to prion diseases. For prion AGAAAAGA amyloid fibril structure, basing on the common property of potential energy minimization of HPs, vdws, and HBs, we will present computational molecular structures of prion AGAAAAGA amyloid fibrils.

## 2.1. Review on Materials and Methods, and Results of MM Models

### 2.1.1. Hybrid Method of Steepest Descent – Conjugate Gradient with Simulated Annealing

X-ray crystallography finds the X-ray final structure of a protein, which usually need refinements using a simulated annealing protocol in order to produce a better structure. Thus, it is very amenable to use simulated annealing (SA) to format the models constructed. Zhang (2011a, 2011d) presents a hybrid method of global search SA with local steepest descent (SD), conjugate gradient (CG) search. The hybrid method is executed with the following three procedures. (1) Firstly the SD method and then the CG method are executed. These two local search methods are traditional optimization methods. The former has nice convergence but is slow when close to minimums. The latter is efficient but its gradient RMS and GMAX gradient (Case et al., 2010) do not have a good convergence. (2) When models cannot be optimized further, we employ standard SA global search procedure. (3) Lastly, the SD and CG methods are used to refine the models. The PDB (Berman et al., 2000) templates used in (Zhang, 2011a, 2011d) are 2OKZ.pdb, 2ONW.pdb, 2OLX.pdb, 2OMQ.pdb, 2ON9.pdb, 2ONV.pdb, 2ONA.pdb, 1XYO.pdb, 2OL9.pdb, 2OMN.pdb, 2ONX.pdb, 2OMP.pdb, 1YJP.pdb of (Sawaya et al., 2007), but only the 2OMP and 1YJP template-based three MM-Models (Fig. 6a~6c in (Zhang, 2011a)) are successfully passed through the SDCG-SA-SDCG computational procedures.

### 2.1.2. Hybrid Method of Discrete Gradient with Simulated Annealing

Zhang et al. (2011a, 2011d) used 3FVA.pdb as the pdb template to build two MM-Models (Figs. 11~12 in (Zhang et al., 2011a)). The Models were built using a hybrid SA Discrete Gradient (DG (Bagirov et al., 2008)) method. Then the Models were optimized using SDCG-SA-SDCG methods as in (Zhang, 2011a).

### 2.1.3. Computational Method of Canonical Dual Global Optimization Theory

Zhang et al. (2011, 2011d) used 3NHC.pdb, 3NVF/G/H/E.pdb templates to build several MM-Models (Figs. 9~11 in (Zhang et al., 2011b), and Figs. 5~8 in (Zhang, 2011b)). These Models were built in the use of canonical dual global optimization theory (Gao et al., 2012; Gao and Wu, 2012; Gao, 2000) and then refined by SDCG-SA-SDCG methods as in (Zhang, 2011a).

## 2.2. New Material and Method, and New MM-Models

### 2.2.1 New Material

This Chapter uses a suitable pdb file template 3NHD.pdb (the GYVLGS segment 127-132 from human prion with V129 (Apostol et al., 2010) from the Protein Data Bank to build MM-Models of AGAAAAGA amyloid fibrils for prions.

### 2.2.2. New Computational Method - Computational Method of Simulated Annealing Evolutionary Computations

The computational methods used to build the new MM-Models will be simulated annealing evolutionary computations (SAEC), where SAECs were got from the hybrid algorithms of (Abbass et al., 2003) by simply replacing the DG method by the SA algorithm of (Bagirov and Zhang, 2003) and numerical computational results show that SAECs can successfully pass the test of more than 40 well-known benchmark global optimization problems (Zhang, 2011c).

### 2.2.3. New MM-Models

The atomic-resolution X-ray structure of 3NHD.pdb is a steric zipper, with strong vdw interactions between β-sheets and HBs to maintain the β-strands (Fig. 2).

Fig. 2: Protein fibril structure of human V129 prion GYVLGS (127–132) (PDB ID: 3NHD). The dashed lines denote the hydrogen bonds. A, B ... K, L denote the chains of the fibril.

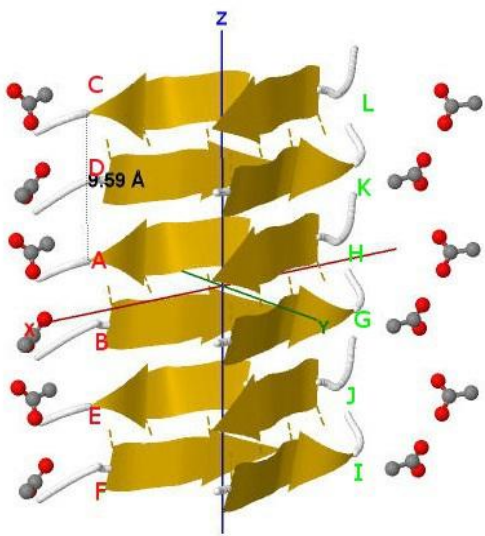

By observations of the 3rd column of coordinates of 3NHD.pdb and Fig. 2, G(H) chains (i.e. β-sheet 2) of 3NHD.pdb can be calculated from A(B) chains (i.e. β-sheet 1) by Eq. 1 and other chains can be calculated by Eqs. 2~3:

$$G(H) = \begin{pmatrix} -1 & 0 & 0 \\ 0 & 1 & 0 \\ 0 & 0 & -1 \end{pmatrix} A(B) + \begin{pmatrix} -20.5865 \\ 9.48 \\ 0.0 \end{pmatrix}, \quad (1)$$

$$K(L) = G(H) + \begin{pmatrix} 0.0 \\ 0.0 \\ 9.59 \end{pmatrix}, I(J) = G(H) + \begin{pmatrix} 0.0 \\ 0.0 \\ -9.59 \end{pmatrix}, \quad (2)$$

$$C(D) = A(B) + \begin{pmatrix} 0.0 \\ 0.0 \\ 9.59 \end{pmatrix}, E(F) = A(B) + \begin{pmatrix} 0.0 \\ 0.0 \\ -9.59 \end{pmatrix}. \quad (3)$$

Basing on the template 3NHD.pdb from the Protein Data Bank, three prion AGAAAAGA palindrome amyloid fibril models - an AGAAAA model (Model 1), a GAAAAG model (Model 2), and an AAAAGA model (Model 3) - will be successfully constructed in this Chapter. Because the template is a segment of 6 residues, the three shorter prion fragments are selected. This Chapter does not perform calculations on the full AGAAAAGA. Chains AB of Models 1~3 were respectively got from AB chains of 3NHD.pdb using the mutate module of the free package Swiss-PdbViewer (SPDBV Version 4.01) (http://spdbv.vital-it.ch). It is pleasant to see that almost all the hydrogen bonds are still kept after the mutations; thus we just need to consider the vdw contacts only. Making mutations for GH chains of 3NHD.pdb, we can get the GH chains of Models 1~3. However, the vdw contacts between Chain A and Chain G, between B chain and H chain are too far at this moment (Figs. 3~5).

Fig. 3: At initial state, the vdw contacts between AB chains (β-sheet 1) and GH chains (β-sheet 2) of Model 1 are very far.

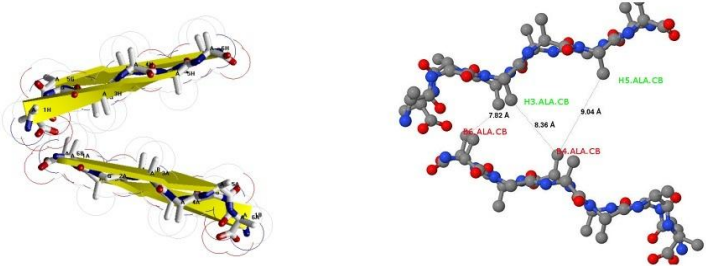

Fig. 4: At initial state, the vdw contacts between AB chains (β-sheet 1) and GH chains (β-sheet 2) of Model 2 are very far.

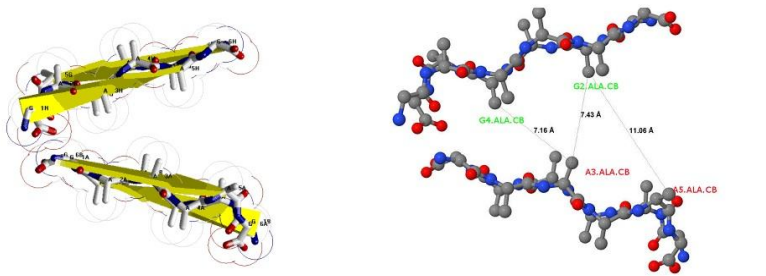

Fig. 5: At initial state, the vdw contacts between AB chains (β-sheet 1) and GH chains (β-sheet 2) of Model 3 are very far.

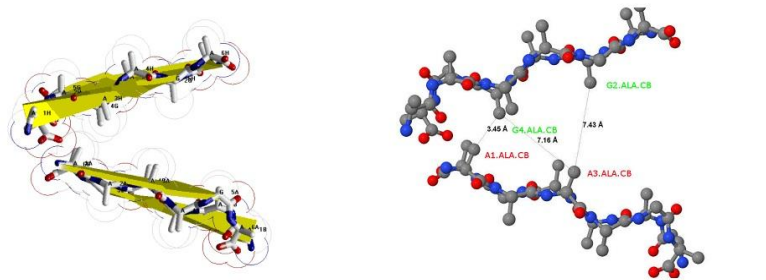

Seeing Figs. 3~5, we may know that for Model 1 at least 3 vdw interactions B6.ALA.CB-H3.ALA.CB-B4.ALA.CB-H5.ALA.CB should be maintained (their distances in Fig. 3 are 7.82, 8.36, 9.04 angstroms respectively), for Model 2 at least 3 vdw interactions G4.ALA.CB-A3.ALA.CB-G2.ALA.CB-A5.ALA.CB should be maintained (their distances in Fig. 4 are 7.16, 7.43, 9.31 angstroms respectively), and for Model 3 at least 3 vdw interactions A1.ALA.CB-G4.ALA.CB-A3.ALA.CB-G2.ALA.CB should be maintained (their distances in Fig. 5 are 3.45, 7.16, 7.43 angstroms respectively). For Model 1, fixing the coordinates of B6.ALA.CB and B4.ALA.CB, letting the coordinates of H3.ALA.CB and H5.ALA.CB be variables, we may get a simple Lennard-Jones (LJ) potential energy minimization problem just with 6 variables (see Eq. 9). Similarly, for Model 2 fixing the coordinates of A3.ALA.CB and A5.ALA.CB, letting the coordinates of G4.ALA.CB and G2.ALA.CB be variables, we may get a simple LJ potential energy minimization problem just with 6 variables (see Eq. 10); for Model 3, fixing the coordinates of A1.ALA.CB and A3.ALA.CB, letting the coordinates of G4.ALA.CB and G2.ALA.CB be variables, we may get a simple LJ potential energy minimization problem with 6 variables (see Eq. 11).

The vdw contacts of atoms are described by the LJ potential energy:

$$V_{LJ}(r) = 4\varepsilon\ [(\sigma/r)^{12} - (\sigma/r)^{6}], \tag{4}$$

where $\varepsilon$ is the depth of the potential well and $\sigma$ is the atom diameter; these parameters can be fitted to reproduce experimental data or deduced from results of accurate quantum chemistry calculations. The $(\sigma/r)^{12}$ term describes repulsion and the $-(\sigma/r)^6$ term describes attraction. If we introduce the coordinates of the atoms whose number is denoted by N and let $\varepsilon=\sigma=1$ be the reduced units, the Eq. 4 becomes into

$$f(x) = 4 \sum_{i=1}^{N} \sum_{j=1, j<i}^{N} (1/t_{ij}^6 - 1/t_{ij}^3), \tag{5}$$

where $t_{ij} = (x_{3i-2} - x_{3j-2})^2 + (x_{3i-1} - x_{3j-1})^2 + (x_{3i} - x_{3j})^2$, $(x_{3i-2}, x_{3i-1}, x_{3i})$ is the coordinates of atom i, $N \geq 2$. The minimization of LJ potential $f(x)$ on $R^n$ (where n = 3N) is an optimization problem:

$$\min f(x) \text{ subject to } x \in R^{3N}. \tag{6}$$

Similarly as Eq. 4, i.e. the potential energy for the vdw interactions between β-sheets:

$$V_{LJ}(r) = A/r^{12} - B/r^6, \tag{7}$$

the potential energy for the HBs between the β-strands has the formula

$$V_{HB}(r) = C/r^{12} - D/r^{10}, \tag{8}$$

where A, B, C, D are given constants. Thus, the amyloid fibril molecular modeling problem can be deduced into the problem to solve the mathematical optimization problem Eq. 6. Seeing Fig. 6, we may know that the optimization problem Eq. 6 reaches its optimal value at the bottom of the LJ potential well, where the distance between two atoms equals to the sum of vdw radii of the atoms. In this Chapter, the sum of the

Fig. 6: The Lennard-Jone Potential (Eqs.4 and 7) (This Fig. can be found in website homepage.mac.com/swain/CMC/DDResources/mol_interactions/molecular_interactions.html).

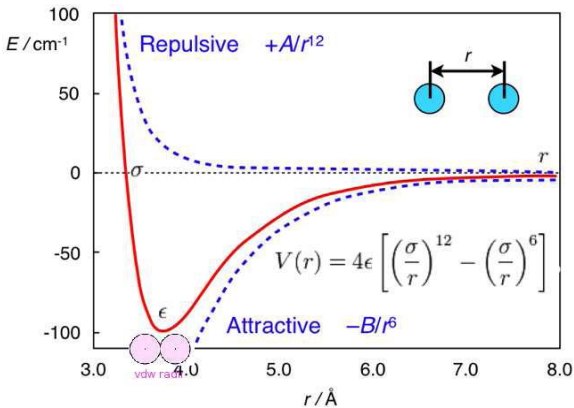

vdw radii is the twice of the vdw radius of Carbon atom, i.e. 3.4 angstroms. The optimization problem Eq. 6 is a nonconvex complex optimization problem. By the observation from Fig. 6, we may solve its simple but equal convex-and-smooth least square optimization problem (or the so-called distance geometry problem or sensor network problem) with a slight perturbation if data for three atoms violate the triangle inequality. The following three optimization problems for Models 1~3 respectively are:

min $f(x)= \frac{1}{2}\{(x_{11}+16.359)^2 + (x_{12}-9.934)^2 + (x_{13}+3.526)^2 - 3.4^2\}^2 + \frac{1}{2}\{(x_{21}+9.726)^2 + (x_{22}-8.530)^2 + (x_{23}+3.613)^2 - 3.4^2\}^2 + \frac{1}{2}\{(x_{11}+9.726)^2 + (x_{12}-8.530)^2 + (x_{13}+3.613)^2 - 3.4^2\}^2$ with initial solution (-12.928, 12.454, 3.034; -6.635, 14.301, 2.628),  (9)

min $f(x)= \frac{1}{2}\{(x_{11}+8.655)^2 + (x_{12}-8.153)^2 + (x_{13}-1.770)^2 - 3.4^2\}^2 + \frac{1}{2}\{(x_{21}+8.655)^2 + (x_{22}-8.153)^2 + (x_{23}-1.770)^2 - 3.4^2\}^2 + \frac{1}{2}\{(x_{21}+2.257)^2 + (x_{22}-6.095)^2 + (x_{23}-3.078)^2 - 3.4^2\}^2$ with initial solution (-13.909, 12.227, -0.889; -7.439, 14.419, -2.033),  (10)

min $f(x)= \frac{1}{2}\{(x_{11}+15.632)^2 + (x_{12}-9.694)^2 + (x_{13}-0.687)^2 - 3.4^2\}^2 + \frac{1}{2}\{(x_{11}+8.655)^2 + (x_{12}-8.153)^2 + (x_{13}-1.770)^2 - 3.4^2\}^2 + \frac{1}{2}\{(x_{21}+8.655)^2 + (x_{22}-8.153)^2 + (x_{23}-1.770)^2 - 3.4^2\}^2$ with initial solution (-13.909, 12.227, -0.889; -7.439, 14.419, -2.033).  (11)

We may use any optimization algorithms or packages to easily solve problems Eqs. 9~11 and get their respective global optimal solutions (-13.062, 9.126, -3.336; -12.344, 6.695, -2.457), (-11.275, 6.606, 3.288; -5.461, 7.124, 2.424), (-12.149, 8.924, 1.229; -9.256, 11.007, 3.517), which were got by the SAEC algorithms in this Chapter. Input these global optimal solutions into Eq. 1, take average and tests then we get Eq. 12:

$$G(H) = \begin{pmatrix} -1 & 0 & 0 \\ 0 & 1 & 0 \\ 0 & 0 & -1 \end{pmatrix} A(B) + \begin{pmatrix} -20.2788 \\ -0.0821 \\ 0.5609 \end{pmatrix}. \qquad (12)$$

By Eq. 12, we can get close vdw contacts in Figs. 7~9.

*Fig. 7: After LJ potential energy minimization, the vdw contacts of Model 1 become very closer (the distances are illuminated by the overlap of border of CB atoms' surface).*

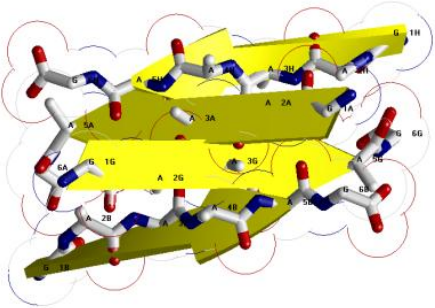

*Fig. 8: After LJ potential energy minimization, the vdw contacts of Model 2 become very closer (the distances are illuminated by the overlap of border of CB atoms' surface).*

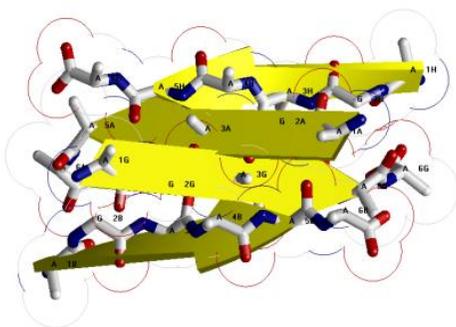

*Fig. 9: After LJ potential energy minimization, the vdw contacts of Model 3 become very closer (the distances are illuminated by the overlap of border of CB atoms' surface).*

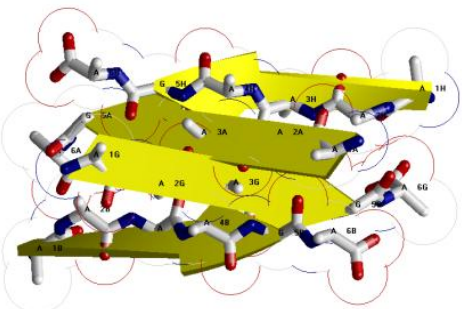

From Figs. 3~5 to Figs. 7~9, we may see that the Optimization algorithm works and the computational experiences show us we had better at least define two sensors and two anchors in order to form a zipper between the two β-sheets. Next, in order to remove very close bad contacts, we relax Figs. 7~9 by a slight SDCG-Optimization in the use of Amber 11 (Case et al., 2010) and we get the optimized MM-Models 1~3. The other CDEF and LKJI chains can be got by parallelizing ABGH chains in the use of mathematical Eqs. 2~3. The new amyloid fibril models are useful for the drive to find treatments for prion diseases in the field of medicinal chemistry. The computational algorithms presented in this Chapter and their references therein are useful in materials science, drug design, etc.

Because Eqs. 9~11 are optimization problems with 6 variables only and these optimization problems are to minimize fourth-order polynomials, the proposed SAEC method and other computational methods can easily get the same optimal solutions to optimize the above three models.

### 2.2.4. The Practical LBFGS Quasi-Newtonian Method

Energy minimization (EM), with the images at the endpoints fixed in space, of the total system energy provides a minimum energy path. EM can be done using SD, CG, and LBFGS (Limited-memory Broyden-Fletcher-Goldfarb-Shanno). SD is robust and easy to implement but it is not most efficient especially when closer to minimum. CG is slower than SD in the early stages but more efficient when closer to minimum. The hybrid of SD-CG will make SD more efficient than SD or CG alone. However, CG cannot be used to find the minimization energy path, for example, when "forces are truncated according to the tangent direction, making it impossible to define a Lagrangian" (Chu et al., 2003). In this case, the powerful and faster quasi-Newtonian method (e.g. the LBFGS quasi-Newtonian minimizer) can be used (Chu et al., 2003; Liu and Nocedal, 1989; Nocedal and Morales, 2000; Byrd et al., 1995; Zhu et al., 1997). We briefly introduce the LBFGS quasi-Newtonian method as follows.

Newton's method in optimization explicitly calculates the Hessian matrix of the second-order derivatives of the objective function and the reverse of the Hessian matrix (Dennis et al., 1996). The convergence of this method is quadratic, so it is faster than SD or CG. In high dimensions, finding the inverse of the Hessian is very expensive. In some cases, the Hessian is a non-invertible matrix, and furthermore in some cases, the Hessian is symmetric indefinite. Qusi-Newton methods thus appear to overcome all these shortcomings.

Quasi-Newton methods (a special case of variable metric methods) are to approximate the Hessian. Currently, the most common quasi-Newton algorithms are the SR1 formula, the BHHH method, the widespread BFGS method and its limited /low-memory extension LBFGS, and Broyden's methods (http://en.wikipedia.org/wiki/Quasi-Newton_method). In Amber (Case et al., 2010) and Gromacs (van der Spoel et al., 2010), LBFGS is used, and the hybrid of LBFGS with CG - a Truncated Newton linear CG method with optional LBFGS Preconditioning (Nocedal and Morales, 2000) - is used in Amber (Case et al., 2010).

### 2.3. New Thinking about the Construction of 3D-Structure of a Protein

If a NMR or X-ray structure of a protein has not been determined and stored in PDB bank yet, we still can easily get the 3D-structural frame of the protein. For example, before 2005 when we did not know the NMR structure of rabbit prion protein, we could get its homology model structure using the NMR structure of the human prion protein (PDB id: 1QLX) as the template (Zhang et al., 2006). We may use the homology structure to determine the 3D-structural frame of a protein when its NMR or X-ray structure has not been determined yet. The determination is an optimization problem described as follows.

"Very often in a structural analysis, we want to approximate a secondary structural element with a single straight line" (Burkowski, 2009: page 212). For example, Fig. 10 uses two straight lines that act as the longitudinal axis of β-Strand A (i.e. A chain), β-Strand B (i.e. B chain) respectively. Each straight line should be positioned among the Cα atoms so that it is closest to all these Cα atoms in a least-squares sense, which is to minimize the sum of the squares of the perpendicular

Fig. 10: The 3D-structural frame of AB chains of Model 1 in Fig. 3 with two β-Strands.

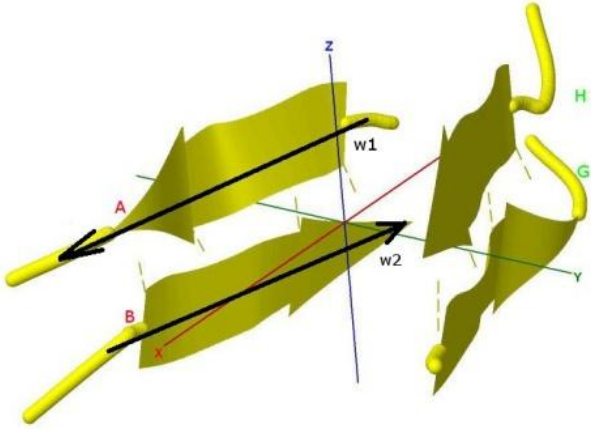

distances ($d_i$) from the Cα atoms to the strand/helix axis:

$$S^* = \min\ S = \sum_{i=1}^{N} \|d_i\|^2. \tag{13}$$

Define the vector $w=(w_x, w_y, w_z)^T$ for the axis. Then $d_i$ represents the perpendicular vector going from Cα atom $a$ to the axis:

$$\|d_i\|^2 = \|a^{(i)}\|^2 \sin^2\theta_i = \|a^{(i)}\|^2 (1-\cos^2\theta_i) = \|a^{(i)}\|^2 \{1 - (a^{(i)T} w)^2 / (\|a^{(i)}\|^2 \|w\|^2)\} = \{a^{(i)}_x{}^2 + a^{(i)}_y{}^2 + a^{(i)}_z{}^2\} \{1 - (a^{(i)}_x w_x + a^{(i)}_y w_y + a^{(i)}_z w_z)^2 / [(a^{(i)}_x{}^2 + a^{(i)}_y{}^2 + a^{(i)}_z{}^2)(w_x{}^2 + w_y{}^2 + w_z{}^2)]\}. \tag{14}$$

According to Eqs. 13~14, for the β-Strand A – β-Strand B of AB chains, we get the following two optimization problems for Model 1 respectively:

$\min S_A = ((-16.196)^2 + 8.315^2 + 1.061^2)\{1- (-16.196w_x + 8.315w_y + 1.061w_z)^2 /[((-16.196)^2 + 8.315^2 + 1.061^2)(w_x{}^2 + w_y{}^2 + w_z{}^2)]\} +$

$((-12.977)^2 + 6.460^2 + 1.908^2)\{1- (-12.977w_x + 6.460w_y + 1.908w_z)^2 /[((-12.977)^2 + 6.460^2 + 1.908^2)(w_x{}^2 + w_y{}^2 + w_z{}^2)]\} +$

$((-9.178)^2 + 6.745^2 + 1.448^2)\{1- (-9.178w_x + 6.745w_y + 1.448w_z)^2 /[((-9.178)^2 + 6.745^2 + 1.448^2)(w_x{}^2 + w_y{}^2 + w_z{}^2)]\} +$

$((-6.455)^2 + 4.112^2 + 1.558^2)\{1- (-6.455w_x + 4.112w_y + 1.558w_z)^2 /[((-6.455)^2 + 4.112^2 + 1.558^2)(w_x{}^2 + w_y{}^2 + w_z{}^2)]\} +$

$((-3.006)^2 + 5.750^2 + 1.782^2)\{1- (-3.006w_x + 5.750w_y + 1.782w_z)^2 /[((-3.006)^2 + 5.750^2 + 1.782^2)(w_x{}^2 + w_y{}^2 + w_z{}^2)]\} +$

$((-1.226)^2 + 2.750^2 + 0.233^2)\{1- (-1.226w_x + 2.750w_y + 0.233w_z)^2 /[((-1.226)^2 + 2.750^2 + 0.233^2)(w_x{}^2 + w_y{}^2 + w_z{}^2)]\},$ (15)

$\min S_B = ((-0.959)^2 + 2.950^2 +(-4.817)^2)\{1- (-0.959w_x + 2.950w_y - 4.817w_z)^2 /[((-0.959)^2 + 2.950^2 +(-4.817)^2)(w_x{}^2 + w_y{}^2 + w_z{}^2)]\} +$

$((-3.465)^2 + 4.999^2 +(-2.846)^2)\{1- (-3.465w_x + 4.999w_y - 2.846w_z)^2 /[((-3.465)^2 + 4.999^2 + (-2.846)^2)(w_x{}^2 + w_y{}^2 + w_z{}^2)]\} +$

$( (-7.213)^2 + 4.412^2 +(-3.340)^2 ) \{ 1- (-7.213w_x + 4.412w_y -3.340w_z )^2 /[( (-7.213)^2 + 4.412^2 + (-3.340)^2 ) (w_x^2 + w_y^2 + w_z^2)]\} +$

$( (-9.954)^2 + 7.078^2 +(-3.168)^2 ) \{ 1- (-9.954w_x + 7.078w_y -3.168w_z )^2 /[( (-9.954)^2 + 7.078^2 + (-3.168)^2 ) (w_x^2 + w_y^2 + w_z^2)]\} +$

$( (-13.660)^2 + 6.241^2 +(-3.137)^2 ) \{ 1- (-13.660w_x + 6.241w_y -3.137w_z )^2 /[( (-13.660)^2 + 6.241^2 + (-3.137)^2 ) (w_x^2 + w_y^2 + w_z^2)]\} +$

$( (-16.702)^2 + 8.507^2 +(-3.074)^2 ) \{ 1- (-16.702w_x + 8.507w_y -3.074w_z )^2 /[( (-16.702)^2 + 8.507^2 + (-3.074)^2 ) (w_x^2 + w_y^2 + w_z^2)]\}.$  (16)

We solve Eqs. 15~16 (taking the average of the coordinates of Cα atoms as initial solutions), getting their optimal solutions **w1** = (-10.751, 6.428, 1.411)$^T$, **w2** = (-7.960, 4.579, -2.256)$^T$ respectively (Fig. 10). We may use the vectors **w1, w2** and Eq. 12 to construct Chains GH and then build an optimal Model 1 (Aqvist, 1986; Abagyan and Maiorov, 1988; Orengo et al., 1992; Young et al., 1999; Foote and Raman, 2000). In (Burkowski, 2009: pages 213-216), $w_x^2 + w_y^2 + w_z^2 = 1$ (i.e. **w** is a unit vector) is restrained and Eq. 13 becomes into a problem to seek the smallest eigenvalue (S*) and its corresponding eigenvector **w** of the following matrix:

$$\begin{pmatrix} \sum_{i=1}^{N}(a_y^{(i)})^2 + (a_z^{(i)})^2 & -\sum_{i=1}^{N} a_x^{(i)} a_y^{(i)} & -\sum_{i=1}^{N} a_z^{(i)} a_x^{(i)} \\ -\sum_{i=1}^{N} a_x^{(i)} a_y^{(i)} & \sum_{i=1}^{N}(a_z^{(i)})^2 + (a_x^{(i)})^2 & -\sum_{i=1}^{N} a_y^{(i)} a_z^{(i)} \\ -\sum_{i=1}^{N} a_z^{(i)} a_x^{(i)} & -\sum_{i=1}^{N} a_y^{(i)} a_z^{(i)} & \sum_{i=1}^{N}(a_x^{(i)})^2 + (a_y^{(i)})^2 \end{pmatrix}.$$

This matrix is symmetric and positive definite, and its eigenvectors form an orthogonal basis for the set of atoms under consideration. In physics, it is called the inertial tensor involving studies of rotational inertia and its eigenvectors are called the principle axes of inertia. Furthermore, we may also notice that Eq. 13 can be rewritten as

min $(\sum_{i=1}^{N} \|d_i\|^2)^2$    subject to    $w^T w = 1$,  (17)

where $\|d_i\|^2 = (a^{(i)2}_x + a^{(i)2}_y + a^{(i)2}_z)(w_x^2 + w_y^2 + w_z^2) - (a^{(i)}_x w_x + a^{(i)}_y w_y + a^{(i)}_z w_z)^2$. Thus, Eq. 17 can be easily solved by the canonical dual global optimization theory (Gao et al., 2012; Gao and Wu, 2012; Gao, 2000), by the ways of solving the canonical dual of Eq. 17 or solving the quadratic differential equations of the prime-dual Gao-Strang complementary function (Gao et al., 2012; Gao and Wu, 2012; Gao, 2000) through some ordinary or partial differential equation computational strategies.

## 3. CONCLUSIONS

To date the hydrophobic region AGAAAAGA palindrome (113-120) of the unstructured N-terminal region (1-123) of prions has little existing experimental structural data available. This Chapter successfully constructs three molecular structure models for AGAAAAGA palindrome (113-120) by using some suitable template 3NHD.pdb from Protein Data Bank and refinement of the Models with several optimization techniques within AMBER 11. These models should be very helpful for the experimental studies of the hydrophobic region AGAAAAGA palindrome of prion proteins (113-120) when the NMR or X-ray molecular structure of prion AGAAAAGA peptide has not been easily determined yet. These constructed Models for amyloid fibrils may be useful for the goals of medicinal chemistry.

This Chapter also introduces numerous practical computational approaches to construct the molecular models when it is difficult to obtain atomic-resolution structures of proteins with traditional experimental methods of X-ray and NMR etc, due to the unstable, noncrystalline and insoluble nature of these proteins. Known structures can be perfectly reproduced by these computational methods, which can be compared with contemporary methods. As we all know, X-ray crystallography finds the X-ray final structure of a protein, which usually need refinements using a SA protocol in order to produce a better structure. SA is a global search procedure and usually it is better to hybrid with local search procedures. Thus, the computational methods introduced in this Chapter should be better than SA along to refine X-ray final structures.

ACKNOWLEDGEMENTS

This research was supported by a Victorian Life Sciences Computation Initiative (VLSCI) grant numbers VR0063 & 488 on its Peak Computing Facility at the University of Melbourne, an initiative of the Victorian Government. The author appreciates the editors for their comments.

REFERENCES

Abagyan R.A., Maiorov V.N. (1988). "A simple quantitative representation of polypeptide chain folds: Comparison of protein tertiary structures". *J. Biomol. Struct. Dyn.* Vol 5, pp 1267-1279.

Abbass H.A., Bagirov A.M., Zhang J.P. (2003). "The discrete gradient evolutionary strategy method for global optimization". *The IEEE Congress on Evolutionary Computation (CEC2003, Australia)*, IEEE-Press, Vol 1, pp 435-442.

Apostol M.I., Sawaya M.R., Cascio D., Eisenberg D. (2010). "Crystallographic studies of prion protein (PrP) segments suggest how structural changes encoded by polymorphism at residue 129 modulate susceptibility to human prion disease". *J. Biol. Chem.* Vol 285, pp 29671-29675.

Aqvist J. (1986). "A simple way to calculate the axes of an α-helix". *Computers & Chemistry* Vol 10, pp 97-99.

Bagirov A.M., Karasozen B., Sezer M. (2008). "Discrete gradient method: a derivative free method for nonsmooth optimization". *J. Opt. Theor. Appl.* Vol 137, pp 317-334.

Bagirov A.M., Zhang J.P. (2003). "Comparative analysis of the cutting angle and simulated annealing methods in global optimization". *Optimization* Vol 52, No 4-5, pp 363-378.

Berman H.M., Westbrook J., Feng Z., Gilliland G., Bhat T. N., Weissig H., Shindyalov I. N., Bourne P. E. (2000). "The protein data bank". *Nucleic Acids Res.* Vol 28, pp 253-242.

Brown D.R. (2000). "Prion protein peptides: optimal toxicity and peptide blockade of toxicity". *Mol. Cell. Neurosci.* Vol 15, pp 66-78.

Brown D.R. (2001). 'Microglia and prion disease". *Microsc. Res. Tech.* Vol 54, pp 71-80.

Brown D.R., Herms J., Kretzschmar H.A. (1994). "Mouse cortical cells lacking cellular PrP survive in culture with a neurotoxic PrP fragment". *Neuroreport* Vol 5, pp 2057-2060.

Byrd R.H., Lu P., Nocedal J. (1995). "A limited memory algorithm for bound constrained optimization". *SIAM J. Scientif. Statistic. Comput.* Vol 16, pp 1190-1208.

Burkowski F.J. (2009). *Structural Bioinformatics: An Algorithmic Approach*, CRC Press, ISBN 9781584886839 (Hardcover).

Cappai R., Collins S.J. (2004). "Structural biology of prions", in *Prions – A Challenge for Science, Medicine and the Public Health System* (Rabenau H.F., Cinatl J., Doerr H.W. (eds)), Basel, Karger, Vol 11, pp 14-32.

Case, D.A., Darden, T.A., Cheatham, T.E., Simmerling, III C.L., Wang, J., Duke, R.E., Luo, R., Walker, R.C., Zhang, W., Merz, K.M., Roberts, B.P., Wang, B., Hayik, S., Roitberg, A., Seabra, G., Kolossvary, I., Wong, K.F., Paesani, F., Vanicek, J., Liu, J., Wu, X., Brozell, S.R., Steinbrecher, T., Gohlke, H., Cai, Q., Ye, X., Wang, J., Hsieh, M.-J., Cui, G., Roe, D.R., Mathews, D.H., Seetin, M.G., Sagui, C., Babin, V., Luchko, T., Gusarov, S., Kovalenko, A., Kollman, P.A. (2010). *AMBER 11*, University of California, San Francisco.

Chabry J., Caughey B., Chesebro B. (1998). "Specific inhibition of in vitro formation of protease-resistant prion protein by synthetic peptides". *J Biol. Chem.* Vol 273 No 21, pp 13203-13207.

Chan J.C.C.(2011). "Steric zipper formed by hydrophobic peptide fragment of Syrian hamster prion protein". *Biochem.* Vol 50, No 32, pp:6815-6823.

Cheng H.M., Tsai T.W.T., Huang W.Y.C., Lee H.K., Lian H.Y., Chou F.C., Mou Y., Chu J., Trout B.L., Brooks B.R. (2003). "A super-linear minimization scheme for the nudged elastic band method". *J. Chem. Phys.* Vol 119, pp 12708-12717.

Daude N. (2004). "Prion diseases and the spleen". *Viral Immunol.* Vol 17, pp 334-349.

Dennis J.E., Robert J.R. Schnabel B. (1996). *Numerical Methods for Unconstrained Optimization and Nonlinear Equations*, SIAM.

Foote J., Raman A. (2000). "A relation between the principal axes of inertia and ligand binding". *Proc. Natl. Acad. Sci. USA* Vol 97, pp 978-983.

Gao D.Y. (2000). Duality *Principles in Nonconvex Systems: Theory, Methods and Applications*, Kluwer Academic Publishers, Dordrecht, Boston, London, ISBN 9780792361459.

Gao D.Y., Ruan N., Pardalos P.M. (2012). "Canonical dual solutions to sum of fourth-order polynomials minimization problems with applications to sensor network localization", in *Sensors: Theory, Algorithms, and Applications* (editors: Vladimir L. Boginski, Clayton W. Commander, Panos M. Pardalos and Yinyu Ye), Springer Optimization and Its Applications Vol. 61, ISBN 9780387886183, pp 37-54.

Gao D.Y., Wu C.Z. (2012). "Triality theory for general unconstrained global optimization problems". *J. Glob. Optim.* arXiv:1104.2970v2.

Gasset M., Baldwin M.A., Lloyd D.H., Gabriel J.M., Holtzman D.M., Cohen F., Fletterick R., Prusiner S.B. (1992). "Predicted alpha-helical regions of the prion protein when synthesized as peptides form amyloid". *Proc. Natl. Acad. Sci. USA* Vol 89 No 22, pp 10940-10944.

Gazit E. (2002). "A possible role for pi-stacking in the self-assembly of amyloid fibrils". *FASEB J.* Vol 16 No 1, pp 77-83.

Gilead S., Gazit E. (2004). "Inhibition of amyloid fibril formation by peptide analogues modified with alpha-aminoisobutyric acid". *Angew. Chem. Int. Ed. Engl.* Vol 43 No 31, pp 4041-4044.

Griffith J.S. (1967). "Self-replication and scrapie". *Nature* Vol 215, pp 1043-1044.

Haigh C.L., Edwards K., Brown D.R. (2005). "Copper binding is the governing determinant of prion protein turnover". *Mol. Cell. Neurosci.* Vol 30 No 2, pp 186-196.

Harrison C.F., Lawson V.A., Coleman B.M., Kim Y.S., Masters C.L., Cappai R., Barnham K.J., Hill A.F. (2010). "Conservation of a glycine-rich region in the prion protein is required for uptake of prion infectivity". *J. Biol. Chem.* Vol 285 No 26, pp 20213-20223.

Holscher C., Delius H., Burkle A. (1998). "Overexpression of nonconvertible PrP$^C$ delta114-121 in scrapie-infected mouse neuroblastoma cells leads to trans-dominant inhibition of wild-type PrP$^{Sc}$ accumulation". *J. Virol.* Vol 72, pp 1153-1159.

Jobling M.F., Huang X., Stewart L.R., Barnham K.J., Curtain C., Volitakis I., Perugini M., White A.R., Cherny R.A., Masters C.L., Barrow C.J., Collins S.J., Bush A.I., Cappai R. (2001). "Copper and Zinc binding modulates the aggregation and neurotoxic properties of the prion peptide PrP 106-126". *Biochem.* Vol 40, pp 8073-8084.

Jobling M.F., Stewart L.R., White A.R., McLean C., Friedhuber A., Maher F., Beyreuther K., Masters C.L., Barrow C.J., Collins S.J., Cappai R. (1999). "The hydrophobic core sequence modulates the neurotoxic and secondary structure properties of the prion peptide 106-126". *J. Neurochem.* Vol 73, pp1557-1565.

Jones E.M., Wu B., Surewicz K., Nadaud P.S., Helmus J.J., Chen S., Jaroniec C.P., Surewicz W.K. (2011). "Structural polymorphism in amyloids: new insights from studies with Y145Stop prion protein fibrils". *J Biol Chem.* Vol 286 No 49, pp 42777-42784.

Kourie J.I. (2001). "Mechanisms of prion-induced modifications in membrane transport properties: implications for signal transduction and neurotoxicity". *Chem. Biol. Interact.* Vol 138 No 1, pp 1-26.

Kourie J.I., Kenna B.L., Tew D., Jobling M.F., Curtain C.C., Masters C.L., Barnham K.J., Cappai R. (2003). "Copper modulation of ion channels of PrP[106-126] mutant prion peptide fragments". *J. Membr. Biol.* Vol 193 No 1, pp 35-45.

Kuwata K., Matumoto T., Cheng H., Nagayama K., James T.L., Roder H. (2003). "NMR-detected hydrogen exchange and molecular dynamics simulations provide structural insight into fibril formation of prion protein fragment 106-126". *Proc. Natl. Acad. Sci. USA* Vol 100 No 25, pp 14790-14795.

Laganowsky A., Liu C., Sawaya M.R., Whitelegge J.P., Park J., Zhao M., Pensalfini A., Soriaga A.B., Landau M., Teng P.K., Cascio D., Glabe C., Eisenberg D. (2012). "Atomic view of a toxic amyloid small oligomer". *Science* Vol 335, pp1228-1231.

Lee S.W., Mou Y., Lin S.Y., Chou F.C., Tseng W.H., Chen C.H., Lu C.Y., Yu S.S., Chan J.C. (2008). "Steric zipper of the amyloid fibrils formed by residues 109-122 of the Syrian hamster prion protein". *J. Mol Biol.* Vol 378 No 5, pp 1142-1154.

Liu D.C., Nocedal J. (1989). "On the limited memory method for large scale optimization". *Math. Programming B* Vol 45, pp 503-528.

Ma B.Y., Nussinov R. (2002). "Molecular dynamics simulations of alanine rich β-sheet oligomers: insight into amyloid formation". *Protein Science* Vol 11, pp 2335-2350.

Morley J.F., Brignull H.R., Weyers J.J., Morimoto R.I. (2002). "The threshold for polyglutamine-expansion protein aggregation and cellular toxicity is dynamic and influenced by aging in Caenorhabditis elegans". *Proc. Natl. Acad. Sci. USA* Vol 99 No 16, pp 10417-10422.

Nelson R., Sawaya M.R., Balbirnie M., Madsen A., Riekel C., Grothe R., Eisenberg D. (2005). "Structure of the cross-beta spine of amyloid-like fibrils". *Nature* Vol 435 No 7043, pp 773-778.

Nocedal J., Morales J. (2000). "Automatic preconditioning by limited memory quasi-Newton updating". *SIAM J. Opt.* Vol 10, pp 1079-1096.

Norstrom E.M., Mastrianni J.A. (2005). "The AGAAAAGA palindrome in PrP is required to generate a productive $PrP^{Sc}$-$PrP^{C}$ complex that leads to prion propagation". *J. Biol. Chem.* Vol 280, pp 27236-27243.

Ogayar A., Snchez-Prez M. (1998). "Prions: an evolutionary perspective". *Internatl. Microbiol.* Vol 1, pp 183-190.

Orengo C.A., Brown N.P., Taylor W.R. (1992). "Fast structure alignment for protein databank searching". *Proteins: Struct. Funct, and Genetics* Vol 14, pp 139-167.

Pan K.M., Baldwin M., Nguyen J. (1993). "Conversion of α-helices into β-sheets features in the formation of the scrapie prion proteins". *Proc. Natl. Acad. Sci. USA* Vol 90, pp 10962-10966.

Pawar A.P., Dubay K.F. et al. (2005). "Prediction of aggregation-prone and aggregation-susceptible regions in proteins associated with neurodegenerative diseases". *J. Mol. Biol.* Vol 350 No 2, pp 379-392.

Reilly C.E. (2000). "Nonpathogenic prion protein ($PrP^C$) acts as a cellsurface signal transducer". *J. Neurol.* Vol 247, pp 819-820.

Sasaki K., Gaikwad J., Hashiguchi S., Kubota T., Sugimura K., Kremer W., Kalbitzer H.R., Akasaka K. (2008). "Reversible monomer-oligomer transition in human prion protein". *Prion* Vol 2 No 3, pp 118-122.

Sawaya M.R., Sambashivan S., Nelson R., Ivanova M.I., Sievers S.A., Apostol M.I., Thompson M.J., Balbirnie M., Wiltzius J.J., McFarlane H.T., Madsen A., Riekel C., Eisenberg D. (2007). "Atomic structures of amyloid cross-beta spines reveal varied steric zippers". *Nature* Vol 447 No 7143, pp 453-457.

Sunde M., Serpell L.C. et al. (1997). "Common core structure of amyloid fibrils by synchrotron X-ray diffraction". *J. Mol. Biol.* Vol 273 No 3, pp 729-739.

van der Spoel D., Lindahl E., Hess B., van Buuren A.R., Apol E., Meulenhoff P.J., Tieleman D.P., Sijbers A.L.T.M., Feenstra K.A., van Drunen R., Berendsen H.J.C (2010). *Gromacs User Manual version 4.5.4*, www.gromacs.org.

Wagoner V.A. (2010). "Computer simulation studies of self-assembly of fibril forming peptides with an Intermediate resolution protein model". *PhD thesis*, North Carolina State University, Raleigh, North Carolina.

Wagoner V.A., Cheon M., Chang I., Hall C.K. (2011). "Computer simulation study of amyloid fibril formation by palindromic sequences in prion peptides". *Proteins: Struct., Funct., and Bioinf.* Vol 79 No 7, pp 2132-2145.

Wegner C., Romer A., Schmalzbauer R., Lorenz H., Windl O., Kretzschmar H.A. (2002). "Mutant prion protein acquires resistance to protease in mouse neuroblastoma cells". *J. Gen. Virol.* Vol 83, pp 1237-1245.

Wormell R.L. (1954). *New fibres from proteins*. Academic Press, page 106.

Young M.M., Skillman A.G., Kuntz I.D. (1999). "A rapid method for exploring the protein structure universe". *Proteins* Vol 34, pp 317-332.

Zanuy D., Ma B., Nussinov R. (2003). "Short peptide amyloid organization: stabilities and conformations of the islet amyloid peptide NFGAIL". *Biophys. J.* Vol 84 No 3, pp 1884-1894.

Zhang J.P. (2011a). "Optimal molecular structures of prion AGAAAAGA palindrome amyloid fibrils formatted by simulated annealing". *J. Mol. Model.* Vol 17 No 1, pp 173-179.

Zhang J.P. (2011b). "Atomic-resolution structures of prion AGAAAAGA amyloid fibrils". In *Amyloids: composition, functions, and pathology* (editors: Irene P. Halcheck and Nancy R. Vernon), Hauppauge, N.Y.: Nova Science Publishers, 2011, ISBN: 9781621005384 (hardcover), Chapter 10: 177-186.

Zhang J.P. (2011c). *Derivative-free hybrid methods in global optimization and applications - in December 2010*, Academic Publishing, ISBN 9783845435800, Chapter 4.2: 92-95.

Zhang J.P. (2011d). *Practical global optimization computing methods in molecular modelling - for atom-resolution structures of amyloid fibrils*, Academic Publishing, ISBN 978-3-8465-2139-7, 192 pages.

Zhang J.P., Sun J., Wu C.Z. (2011a). "Optimal atomic-resolution structures of prion AGAAAAGA amyloid fibrils". *J. Theor. Biol.* Vol 279 No 1, pp 17-28.

Zhang J.P., Gao D.G., Yearwood J. (2011b). "A novel canonical dual computational approach for prion AGAAAAGA amyloid fibril molecular modeling". *J. Theor. Biol.* Vol 284 No 1, pp 149-157.

Zhang J.P., Varghese J.N., Epa V.C. (2006). "Studies on the conformational stability of the rabbit prion protein". *CSIRO Preventative Health National Research Flagship Science Retreat, Aitken Hill, Melbourne, 12-15 September 2006*, Poster in Excellence.

Zhang Z.Q., Chen H., Lai L.H. (2007). "Identification of amyloid fibril-forming segments based on structure and residue-based statistical potential". *Bioinf.* Vol 23, pp 2218-2225.

Zhu C., Byrd R.H., Nocedal J. (1997). "L-BFGS-B: Algorithm 778: L-BFGS-B, FORTRAN routines for large scale bound constrained optimization". *ACM Trans. Math. Softw.* Vol 23, pp 550-560.